\documentstyle[amssymb,12pt]{article}

\newtheorem{th}{Theorem}[section]
\newtheorem{lem}[th]{Lemma}

\newtheorem{cor}[th]{Corollary}
\newtheorem{defn}[th]{Definition}
\newenvironment{defn-new}{\begin{defn} \em}{\end{defn}}
\newtheorem{rem}[th]{Remark}
\newenvironment{rem-new}{\begin{rem} \em}{\end{rem}}
\newtheorem{ex}[th]{Example}
\newenvironment{ex-new}{\begin{ex} \em}{\end{ex}}
\newtheorem{prob}[th]{Problem}
\newenvironment{prob-new}{\begin{prob} \em}{\end{prob}}

\newenvironment{notation-new}{\begin{rem} \em}{\end{rem}}

\newenvironment{agr-new}{\begin{rem} \em}{\end{rem}}

\makeatletter \@addtoreset{equation}{section} \makeatother
\begin{document}

\bigskip

\begin{center}
{\large {\bf Improved Chen-Ricci inequality}}\medskip

{\large {\bf for curvature-like tensors and its applications}} \bigskip
\bigskip

Mukut Mani Tripathi\bigskip
\end{center}

{\bf Abstract.} We present Chen-Ricci inequality and improved
Chen-Ricci inequality for curvature like tensors. Applying our
improved Chen-Ricci inequality we study Lagrangian and Kaehlerian
slant submanifolds of complex space forms, and $C$-totally real
submanifolds of Sasakian space forms.
\medskip

 \medskip


\section{Introduction\label{sect-intro}}

Since the celebrated theory of J. F. Nash of isometric immersion of a
Riemannian manifold in to a suitable Euclidean space gives very important
and effective motivation to view each Riemannian manifold as a submanifold
in a Euclidean space, the problem of discovering simple sharp relationships
between intrinsic and extrinsic invariants of a Riemannian submanifold
becomes one of the most fundamental problems in submanifold theory. The main
extrinsic invariant is the squared mean curvature and the main intrinsic
invariants include the classical curvature invariants namely the Ricci
curvature and the scalar curvature. There are also many other important
modern intrinsic invariants of (sub)manifolds introduced by B.-Y. Chen (Chen
2000, \cite{Chen-00-handbook}). \medskip

In 1999, B.-Y. Chen \cite{Chen-99-Glasg} proved a basic inequality involving
the Ricci curvature ${\rm Ric}$ and the squared mean curvature $\Vert H\Vert
^{2}$ of submanifolds in a real space form as follows.

\begin{th}
\label{th-Ricci-RSF} {\rm (\cite[Theorem~4]{Chen-99-Glasg})} Let $M$ be an $%
n $-dimensional submanifold of a real space form $R^{m}(c)$. Then the
following statements are true.

\begin{enumerate}
\item[{\bf (a)}] For each unit vector $X \in T_{p}M$, we have
\begin{equation}
\Vert H\Vert ^{2}\geq \frac{4}{n^{2}}\,\left\{ {\rm Ric}(X) - (n-1)c\right\}
.  \label{eq-Ricci-RSF}
\end{equation}

\item[{\bf (b)}] If $H(p)=0$, then a unit vector $X\in T_{p}M$ satisfies the
equality case of $\left( \ref{eq-Ricci-RSF}\right) $ if and only if $X$
belongs to the relative null space ${\cal N}(p)$ given by
\[
{\cal N}(p) = \left\{ X\in T_{p}M : \sigma (X,Y) = 0\; {\rm for\ all}\; Y\in
T_{p}M\right\}.
\]

\item[{\bf (c)}] The equality case of $( \ref{eq-Ricci-RSF})$ holds for all
unit vectors $X$ $\in T_{p}M$ if and only if either $p$ is a geodesic point
or $n=2$ and $p$ is an umbilical point.
\end{enumerate}
\end{th}

The inequality (\ref{eq-Ricci-RSF}) drew attention of several authors and
they established similar inequalities for different kind of submanifolds in
ambient manifolds possessing different kind of structures. The submanifolds
included mainly invariant, anti-invariant and slant submanifolds, while
ambient manifolds included mainly real space forms, complex space forms and
Sasakian space forms. Thus, after putting an extra condition on the Riemann
curvature tensor of the ambient manifold, like its constancy in the case of
real space forms, the constancy of holomorphic sectional curvature in the
case of complex space forms and the constancy of $\varphi $-holomorphic
sectional curvature in the case of Sasakian space forms; one proves the
results similar to that of \cite[Theorem~4]{Chen-99-Glasg} or \cite[Theorem~1%
]{Chen-00-Arch-Math}. \medskip

Motivated by the result of B.-Y. Chen (\cite[Theorem~4]{Chen-99-Glasg}), in
\cite{Hong-Tri-05-IJPAMS-1} and \cite{Hong-MT-05-SUT}, the authors presented
a general theory for a submanifold of Riemannian manifolds by proving a
basic inequality, now called Chen-Ricci inequality \cite{Tri-08-JAdMS},
involving the Ricci curvature and the squared mean curvature of the
submanifold. The goal was achieved by use of the concept of $k$-Ricci
curvature ($2\leq k\leq n$) in an $n$-dimensional Riemannian manifold \cite%
{Chen-99-Glasg}. It can be noted that a $k$-Ricci curvature is a $(k-1)$%
-Ricci curvature in the sense of H. Wu \cite{Wu-87-IUMJ}. In fact, without
assuming any further condition on the Riemann curvature tensor of the
ambient manifold $\widetilde{M}$, we established a Chen-Ricci inequality
involving Ricci curvature and the squared mean curvature for a submanifold $%
M $ of $\widetilde{M}$ as follows.

\begin{th}
\label{th-Chen-Ricci} Let $M$ be an $n$-dimensional submanifold of a
Riemannian manifold. Then, the following statements are true.

\begin{enumerate}
\item[{\bf (a)}] For $X \in T_{p}^{1}M$, it follows that
\begin{equation}
{\rm Ric}(X) \leq \frac{1}{4}\,n^{2}\Vert H\Vert ^{2} + \widetilde{{\rm Ric}}%
_{(T_{p}M)}(X) ,  \label{eq-Ricci}
\end{equation}
where $\widetilde{{\rm Ric}}_{(T_{p}M)}(X) $ is the $n$-Ricci curvature of $%
T_{p}M$ at $X$ $\in T_{p}^{1}M$ with respect to the ambient manifold $%
\widetilde{M}$.

\item[{\bf (b)}] The equality case of $(\ref{eq-Ricci})$ is satisfied by $%
X\in T_{p}^{1}M$ if and only if
\begin{equation}
\left\{
\begin{array}{lll}
\sigma (X,Y) =0, &  & {\rm for\ all\ }Y\in T_{p}M\; {\rm orthogonal\ to}\;
X, \medskip \\
2\sigma (X,X) = nH(p). &  &
\end{array}
\right.  \label{eq-Ricci-0}
\end{equation}
If $H(p)=0$, then $X\in T_{p}^{1}M$ satisfies the equality case of $(\ref%
{eq-Ricci-RSF})$ if and only if $X\in {\cal N}(p)$.

\item[{\bf (c)}] The equality case of $(\ref{eq-Ricci})$ holds for all $X
\in T_{p}^{1}M$ if and only if either $p$ is a geodesic point or $n=2 $ and $%
p$ is an umbilical point.\medskip
\end{enumerate}
\end{th}

Continuing this study of \cite{Hong-Tri-05-IJPAMS-1} and \cite%
{Hong-MT-05-SUT} we also studied Chen-Ricci inequality for submanifolds in
contact metric manifolds and obtained many interesting results (see \cite%
{Hong-Tri-05-IJPAMS-2}, \cite{Hong-Tri-06-IJMSI}, \cite{Tri-08-JAdMS} and
references cited therein). \medskip

In 2005, Oprea \cite{Oprea-05-arXiv} (see also \cite{Oprea-10-MIA}) proved
Chen-Ricci inequality by using optimization techniques applied in the setup
of Riemannian geometry. He also improved Chen-Ricci inequality for
Lagrangian submanifolds of complex space forms. Later, Deng \cite%
{Deng-09-IEJG} proved the improved Chen-Ricci inequality for Lagrangian
submanifolds of complex space forms just by using some crucial algebraic
inequalities and also discussed the equality case. \medskip

However, improved Chen-Ricci inequalities for Kaehlerian slant submanifolds
of complex space forms and $C$-totally real submanifolds of Sasakian space
forms are not known so far. Even improved Chen-Ricci inequalities in these
two cases cannot be obtained directly from the results of Oprea \cite%
{Oprea-05-arXiv} and Deng \cite{Deng-09-IEJG}. \medskip

Under these circumstances it becomes necessary to give a general theory,
which could be applied to obtain Chen-Ricci inequality and improved
Chen-Ricci (in)equality in different situations. Motivated by \cite%
{Bolt-DFV-09-MIA}, we present Chen-Ricci inequality and improved Chen-Ricci
inequality for curvature like tensors (see Theorems~\ref{th-T-Ricci} and \ref%
{th-T-Ricci-imp}). Then we apply our improved Chen-Ricci inequality for
curvature like tensors in study of Kaehlerian slant submanifolds of complex
space forms and $C$-totally real submanifolds of Sasakian space forms. In
the process of this study, we come across several natural problems, which
can be studied in future. \medskip

The paper is organized as follows. In section~\ref%
{sect-Chen-Ricci-inequality}, first we give concepts related with curvature
like tensors. Next, given an $n$-dimensional Riemannian manifold $\left(
M,g\right) $, a Riemannian vector bundle $\left( B,g_{B}\right) $ over $M$,
a $B$-valued symmetric $\left( 1,2\right) $-tensor field $\zeta $ and a
(curvature-like) tensor field $T$ satisfying the algebraic Gauss equation $(%
\ref{eq-Gauss-alg})$, we establish Chen-Ricci inequality (\ref{eq-T-Ric-ineq}%
) involving $T$-Ricci curvature ${\rm Ric}_{T}$ and $\left\Vert {\rm trace}%
\,\zeta \right\Vert $. In section~\ref{sect-Improved-Chen-Ricci-inequality},
we improve Chen-Ricci inequality (\ref{eq-T-Ric-ineq}) under certain
restrictions on $\zeta $ and obtain improved Chen-Ricci inequality (\ref%
{eq-T-Ric-ineq-imp}) (cf. Theorem~\ref{th-T-Ricci-imp}). In section~\ref%
{sect-Lagrangian-submanifolds}, applying our main Theorem~\ref%
{th-T-Ricci-imp} we obtain improved Chen-Ricci inequality for Lagrangian
submanifolds of a complex space form in Theorem~\ref{th-T-Ricci-imp-Lag}
\cite[Theorem~3.1]{Deng-09-IEJG}, which is an improvement of \cite[Theorem~4]%
{Chen-99-Glasg}. It is known that \cite[Example~3.1]{Deng-09-IEJG} the
Whitney $2$-sphere in ${\Bbb C}^{2}$ satisfies the equality case of the
improved Chen-Ricci inequality (\ref{eq-Ricci-imp-Lag}). In Section~\ref%
{sect-Kaehlerian-slant-submanifolds}, we next apply our Theorem~\ref%
{th-T-Ricci-imp} and obtain improved Chen-Ricci inequality (\ref%
{eq-Ricci-imp-Kaeh-slant}) for Kaehlerian slant submanifolds in complex
space forms. This inequality is an improvement of Chen-Ricci inequality for
Kaehlerian slant submanifolds in complex space forms \cite[Inequality~(2.1)
of Theorem~2.1]{Mat-Mih-Taz-03}. We also note that totally umbilical
Lagrangian submanifolds, of dimension $n\geq 2$, in a complex space form
must be totally geodesic \cite[Theorem 1]{Chen-Ogiue-74-MMJ}. As an
improvement of this result, in Theorem~\ref{th-tuL-1}, we prove that if $M$
is a totally umbilical Lagrangian submanifold of a Kaehler manifold then
either $\dim \left( M\right) =1$ or $M$ is totally geodesic. Next, combining
Theorem~\ref{th-tuL-1} with the result that every totally umbilical proper
slant submanifold of a Kaehler manifold is totally geodesic \cite[Theorem~3.1%
]{Sahin-09-RM}, we also conclude that \label{th-tuL-2 copy(1)} each $n$%
-dimensional ($n\geq 2$) totally umbilical non-invariant slant submanifold
of a $2n$-dimensional Kaehler manifold is always totally geodesic (cf.
Theorem~\ref{th-tuL-2}). We also discover that proper slumbilical surfaces
\cite{Chen-02-Osaka} cannot satisfy the equality case of the improved
Chen-Ricci equality (\ref{eq-Ricci-imp-Kaeh-slant}); thus we propose the
definition and classification of $H$-slumbilical surfaces in complex space
forms (cf.~Problem~\ref{prob-H-slumb-1}). Finally, in section~\ref%
{sect-C-totally-real-submanifolds}, applying our Theorem~\ref{th-T-Ricci-imp}%
, we obtain improved Chen-Ricci inequality (\ref{eq-Ricci-imp-C-tot}) for $C$%
-totally real submanifolds of a Sasakian space form (cf. Theorem~\ref%
{th-T-Ricci-imp-C-tot}), which is an improvement of Chen-Ricci inequality
\cite[Inequality~(2.1) of Theorem~2.1]{Mih-02-JAMS}. Like the concept of $H$%
-umbilical Lagrangian submanifolds \cite{Chen-97-Israel}, we propose the
definition and classification of $H$-umbilical $C$-totally real submanifolds
in Sasakian space forms (cf. Problem~\ref{prob-H-umb-C-tot-real}).

\section{Chen-Ricci inequality\label{sect-Chen-Ricci-inequality}}

Let $(M,g)$ be an $n$-dimensional Riemannian manifold. Let $T$ be a
curvature-like tensor field so that it satisfies the following symmetry
properties
\[
T(X,Y,Z,W)=-\,T(Y,X,Z,W),
\]%
\[
T(X,Y,Z,W)=-\,T(Y,X,W,Z),
\]%
\[
T(X,Y,Z,W)+T(X,Z,W,Y)+T(X,W,Y,Z)=0
\]%
for all vector fields $X$, $Y$, $Z$ and $W$ on $M$. For a curvature-like
tensor field $T$, the $T${\em -sectional curvature} associated with a $2$%
-plane section $\Pi _{2}$ spanned by orthonormal vectors $X$ and $Y$ at $%
p\in M$, is given by \cite{Bolt-DFV-09-MIA}
\[
K_{T}(\Pi _{2})=K_{T}(X\wedge Y)=T(X,Y,Y,X).
\]%
Let $\left\{ e_{1},e_{2},\ldots ,e_{n}\right\} $ be any orthonormal basis of
$T_{p}M$. The $T${\em -Ricci tensor} $\,S_{T}$ is defined by
\[
S_{T}(X,Y)=\sum_{j=1}^{n}T\left( e_{j},X,Y,e_{j}\right) ,\qquad X,Y\in
T_{p}M.
\]%
The $T${\em -Ricci curvature} is given by
\[
{\rm Ric}_{T}(X)=S_{T}(X,X),\qquad X\in T_{p}M.
\]%
We denote the set of unit vectors in $T_{p}M$ by $T_{p}^{1}M$; thus
\[
T_{p}^{1}M=\left\{ X\in T_{p}M\ |\ \left\langle X,X\right\rangle =1\right\} .
\]

If $T$ is replaced by the Riemann curvature tensor $R$, then $T$-sectional
curvature $K_{T}$, $T$-Ricci tensor $S_{T}$ and $T$-Ricci curvature ${\rm Ric%
}_{T}$, become sectional curvature $K$, Ricci tensor $S$ and Ricci curvature
${\rm Ric}$, respectively. \medskip

Let $\left( M,g\right) $ be an $n$-dimensional Riemannian manifold and $%
\left( B,g_{B}\right) $ a Riemannian vector bundle over $M$. If $\zeta $ is
a $B$-valued symmetric $\left( 1,2\right) $-tensor field and $T$ a $\left(
0,4\right) $-tensor field on $M$ such that
\begin{equation}
T(X,Y,Z,W)=g_{B}(\zeta (X,W),\zeta (Y,Z))-g_{B}(\zeta (X,Z),\zeta (Y,W))
\label{eq-Gauss-alg}
\end{equation}%
for all vector fields $X$,$Y$,$Z$,$W$ on $M$, then the equation (\ref%
{eq-Gauss-alg}) is said to be an algebraic Gauss equation. If $T$ is a $%
(0,4) $-tensor field on $M$ which satisfies (\ref{eq-Gauss-alg}) then $T$
becomes curvature-like. A typical example of an algebraic Gauss equation is
given for a submanifold $M$ of an Euclidean space, if $B$ is the normal
bundle, $\zeta $ the second fundamental form and $T$ the curvature tensor.
Some nice situations, in which such $T$ and $\zeta $ satisfy an algebraic
Gauss equation exist, are Lagrangian and Kaehlerian slant submanifolds of
complex space forms, $C$-totally real submanifolds of Sasakian space forms,
and positive definite centroaffine hypersurfaces.

\medskip Now, let $\{e_{1},\ldots ,e_{n}\}$ be an orthonormal basis of the
tangent space $T_{p}M$ and $e_{r}$ belong to an orthonormal basis $%
\{e_{n+1},\ldots ,e_{m}\}$ of the Riemannian vector bundle $(B,g_{B})$ over $%
M$ at $p$. We put
\[
\zeta _{ij}^{r}=g_{B}\left( \zeta \left( e_{i},e_{j}\right) ,e_{r}\right)
,\quad \quad \left\Vert \zeta \right\Vert ^{2}=\sum_{i,j=1}^{n}g_{B}\left(
\zeta \left( e_{i},e_{j}\right) ,\zeta \left( e_{i},e_{j}\right) \right) ,
\]
\[
{\rm trace}\,\zeta =\sum_{i=1}^{n}\zeta \left( e_{i},e_{i}\right) ,\quad
\quad \left\Vert {\rm trace}\,\zeta \right\Vert ^{2}=g_{B}({\rm trace}%
\,\zeta ,{\rm trace}\,\zeta ),
\]
\[
{\cal N}_{\zeta }(p)=\left\{ X\in T_{p}M:\zeta \left( X,Y\right) =0\;{\rm %
for\ all}\;Y\in T_{p}M\right\} .
\]

\begin{th}
\label{th-T-Ricci} Let $\left( M,g\right) $ be an $n$-dimensional Riemannian
manifold, $\left( B,g_{B}\right) $ a Riemannian vector bundle over $M$ and $%
\zeta $ a $B$-valued symmetric $\left( 1,2\right) $-tensor field. Let $T$ be
a curvature-like tensor field satisfying the algebraic Gauss equation $(\ref%
{eq-Gauss-alg})$. Then, the following statements are true\/{\rm :}

\begin{enumerate}
\item[{\rm (a)}] For $X$ $\in T_{p}^{1}M$, it follows that
\begin{equation}
{\rm Ric}_{T}(X)\leq \frac{1}{4}\left\Vert {\rm trace}\,\zeta \right\Vert
^{2}.  \label{eq-T-Ric-ineq}
\end{equation}

\item[{\rm (b)}] The equality case of $\left( \ref{eq-T-Ric-ineq}\right) $
is satisfied by $X\in T_{p}^{1}M$ if and only if
\begin{equation}
\left\{
\begin{array}{lll}
\zeta \left( X,Y\right) =0, &  & {\rm for\ all}\;Y\in T_{p}M\ {\rm such\ that%
}\;g\left( X,Y\right) =0,\medskip \\
\zeta \left( X,X\right) =\displaystyle\frac{1}{2}\,{\rm trace}\,\zeta . &  &
\end{array}%
\right.  \label{eq-T-Ric-ineq-1}
\end{equation}

\item[{\rm (c)}] The equality case of the inequality $\left( \ref%
{eq-T-Ric-ineq}\right) $ is true for all $X\in T_{p}^{1}M$ if and only if
either $\zeta =0$\ or $n=2$ and
\begin{equation}
\zeta _{11}^{r}=\zeta _{22}^{r}=\frac{1}{2}\left( \zeta _{11}^{r}+\zeta
_{22}^{r}\right)  \label{eq-T-Ric-ineq-2}
\end{equation}
\end{enumerate}
\end{th}

\noindent {\bf Proof.} First, we note that
\begin{eqnarray}
\Vert \zeta \Vert ^{2} &=&\frac{1}{2}\left\Vert {\rm trace}\,\zeta
\right\Vert ^{2}+\frac{1}{2}\sum_{r=n+1}^{m}(\zeta _{11}^{r}-\zeta
_{22}^{r}-\cdots -\zeta _{nn}^{r})^{2}  \nonumber \\
&&+\ 2\sum_{r=n+1}^{m}\sum_{j=2}^{n}(\zeta
_{1j}^{r})^{2}-2\sum_{r=n+1}^{m}\sum_{2\leq i<j\leq n}(\zeta _{ii}^{r}\zeta
_{jj}^{r}-(\zeta _{ij}^{r})^{2}).  \label{eq-zeta-square}
\end{eqnarray}%
From (\ref{eq-Gauss-alg}), we get
\begin{equation}
(K_{T})_{ij}=\sum_{r=n+1}^{m}\left( \zeta _{ii}^{r}\zeta _{jj}^{r}-(\zeta
_{ij}^{r})^{2}\right) ,  \label{eq-T-Kij}
\end{equation}%
which implies that
\begin{equation}
\tau _{T}(p)=\frac{1}{2}\,\left\Vert {\rm trace}\,\zeta \right\Vert ^{2} -
\frac{1}{2}\,\Vert \zeta \Vert ^{2}.  \label{eq-T-tau-tr-zeta-sq}
\end{equation}%
From (\ref{eq-T-tau-tr-zeta-sq}) and (\ref{eq-zeta-square}) we get
\begin{eqnarray}
\tau _{T}(p) &=&\frac{1}{4}\,\left\Vert {\rm trace}\,\zeta \right\Vert ^{2}
- \frac{1}{4}\sum_{r=n+1}^{m}(\zeta _{11}^{r}-\zeta _{22}^{r}-\cdots -\zeta
_{nn}^{r})^{2}  \nonumber \\
&&-\sum_{r=n+1}^{m}\sum_{j=2}^{n}(\zeta
_{1j}^{r})^{2}+\sum_{r=n+1}^{m}\sum_{2\leq i<j\leq n}(\zeta _{ii}^{r}\zeta
_{jj}^{r}-(\zeta _{ij}^{r})^{2}).  \label{eq-T-Ric-ineq-3}
\end{eqnarray}%
From (\ref{eq-T-Kij}) we also have
\begin{equation}
\sum_{2\leq i<j\leq n}(K_{T})_{ij}=\sum_{r=n+1}^{m}\ \sum_{2\leq i<j\leq
n}\left( \zeta _{ii}^{r}\zeta _{jj}^{r}-(\zeta _{ij}^{r})^{2}\right) .
\label{eq-T-Ric-ineq-4}
\end{equation}%
From (\ref{eq-T-Ric-ineq-3}) and (\ref{eq-T-Ric-ineq-4}), we obtain
\begin{equation}
{\rm Ric}_{T}\left( e_{1}\right) =\frac{1}{4}\,\left\Vert {\rm trace}\,\zeta
\right\Vert ^{2}-\sum_{r=n+1}^{m}\sum_{j=2}^{n}(\zeta _{1j}^{r})^{2}-\frac{1%
}{4}\sum_{r=n+1}^{m}(\zeta _{11}^{r}-\zeta _{22}^{r}-\cdots -\zeta
_{nn}^{r})^{2}.  \label{eq-T-Ric-ineq-5}
\end{equation}%
Since, we can choose $e_{1}=X$ as any unit vector in $T_{p}^{1}M$, therefore
(\ref{eq-T-Ric-ineq-5}) implies (\ref{eq-T-Ric-ineq}).

To prove the statement (b), assuming $X=e_{1}$, from (\ref{eq-T-Ric-ineq-5}%
), the equality in (\ref{eq-T-Ric-ineq}) is valid if and only if
\begin{equation}
\zeta _{12}^{r}=\cdots =\zeta _{1n}^{r}=0\quad {\rm and}\quad \zeta
_{11}^{r}=\zeta _{22}^{r}+\cdots +\zeta _{nn}^{r},\quad r\in \left\{
n+1,\ldots ,m\right\} ,  \label{eq-T-Ric-ineq-6}
\end{equation}%
which is equivalent to (\ref{eq-T-Ric-ineq-1}).

Now we prove the statement (c). Assuming the equality case of (\ref%
{eq-T-Ric-ineq}) for all unit vectors $X\in T_{p}^{1}M$, in view of (\ref%
{eq-T-Ric-ineq-6}), for each $r\in \left\{ n+1,\ldots ,m\right\} $ it
follows that
\begin{equation}
\zeta _{ij}^{r}=0,\qquad i\neq j.  \label{eq-T-Ric-ineq-7}
\end{equation}%
\begin{equation}
2\zeta _{ii}^{r}=\zeta _{11}^{r}+\zeta _{22}^{r}+\cdots +\zeta
_{nn}^{r},\qquad i\in \left\{ 1,\ldots ,n\right\} ,  \label{eq-T-Ric-ineq-8}
\end{equation}%
From (\ref{eq-T-Ric-ineq-8}), we have
\[
2\zeta _{11}^{r}=2\zeta _{22}^{r}=\cdots =2\zeta _{nn}^{r}=\zeta
_{11}^{r}+\zeta _{22}^{r}+\cdots +\zeta _{nn}^{r}\,,
\]%
which implies that
\[
(n-2)\left( \zeta _{11}^{r}+\zeta _{22}^{r}+\cdots +\zeta _{nn}^{r}\right)
=0.
\]%
Thus, either $\zeta _{11}^{r}+\zeta _{22}^{r}+\cdots +\zeta _{nn}^{r}=0$ or $%
n=2$. If $\zeta _{11}^{r}+\zeta _{22}^{r}+\cdots +\zeta _{nn}^{r}=0,$ then
in view of (\ref{eq-T-Ric-ineq-8}), we get
\[
\zeta _{ii}^{r}=0,\qquad i\in \left\{ 1,\ldots ,n\right\} .
\]%
This together with (\ref{eq-T-Ric-ineq-7}) gives $\zeta _{ij}^{r}=0$ for all
$i,j\in \left\{ 1,\ldots ,n\right\} $ and $r\in \left\{ n+1,\ldots
,m\right\} $, that is, $\zeta =0$. If $n=2$, then from (\ref{eq-T-Ric-ineq-8}%
) we get (\ref{eq-T-Ric-ineq-2}). The proof of the converse part is
straightforward. $\blacksquare $ \medskip

We immediately have the following

\begin{cor}
\label{cor-T-Ricci} Let $\left( M,g\right) $ be an $n$-dimensional
Riemannian manifold, $\left( B,g_{B}\right) $ a Riemannian vector bundle
over $M$ and $\zeta $ a $B$-valued symmetric $\left( 1,2\right) $-tensor
field. Let $T$ be a curvature-like tensor field defined by $(\ref%
{eq-Gauss-alg})$. Then for $X\in T_{p}^{1}M$ any two of the following three
statements imply the remaining one.

\begin{enumerate}
\item[{\rm (a)}] $X$ satisfies the equality case of $(\ref{eq-T-Ric-ineq})$.

\item[{\rm (b)}] ${\rm trace}\,\zeta (p)=0$.

\item[{\rm (c)}] $X\in {\cal N}_{\zeta }(p)$.
\end{enumerate}
\end{cor}

Several results in form of Chen-Ricci inequalities in \cite{Hong-MT-05-SUT},
\cite{Hong-Tri-05-IJPAMS-1}, \cite{Hong-Tri-05-IJPAMS-2}, \cite%
{Hong-Tri-06-IJMSI} and papers cited therein are among others, which can be
proved by suitable applications of Theorem~\ref{th-T-Ricci}. However, in
this paper we are concerned with improved Chen-Ricci inequalities in
different situations. Now, in the following section, we improve the
Chen-Ricci inequality for curvature like tensors satisfying an algebraic
Gauss equation $(\ref{eq-Gauss-alg})$ under certain restrictions on the
tensor field $\zeta$.

\section{Improved Chen-Ricci inequality\label%
{sect-Improved-Chen-Ricci-inequality}}

First, we state following two Lemmas for later use.

\begin{lem}
\label{lem-Deng-1} {\rm (\cite[Lemma 2.2]{Deng-09-IEJG})} Let $f_{1}:{\Bbb R}%
^{n}\rightarrow {\Bbb R}$ be a function defined by
\[
f_{1}\left( a^{1},\ldots ,a^{n}\right)
=a^{1}\sum_{j=2}^{n}a^{j}-\sum_{j=2}^{n}\left( a^{j}\right) ^{2}.
\]%
If $a^{1}+\cdots +a^{n}=2na$, we have
\[
f_{1}\left( a^{1},\ldots ,a^{n}\right) \leq \frac{n-1}{4n}\,\left(
a^{1}+\cdots +a^{n}\right) ^{2}.
\]%
The equality sign holds if and only if
\[
\frac{1}{n+1}\,a^{1}=a^{2}=\cdots =a^{n}=a.
\]
\end{lem}

\begin{lem}
\label{lem-Deng-2} {\rm (\cite[Lemma 2.3]{Deng-09-IEJG})} Let $f_{2}:{\Bbb R}%
^{n}\rightarrow {\Bbb R}$ be a function defined by
\[
f_{2}\left( a^{1},\ldots ,a^{n}\right) =a^{1}\sum_{j=2}^{n}a^{j}-\left(
a^{1}\right) ^{2}.
\]%
If $a^{1}+\cdots +a^{n}=4a$, we have
\[
f_{2}\left( a^{1},\ldots ,a^{n}\right) \leq \frac{1}{8}\,\left( a^{1}+\cdots
+a^{n}\right) ^{2}.
\]%
The equality sign holds if and only if
\[
a^{1}=a,\quad a^{2}+\cdots +a^{n}=3a.
\]
\end{lem}

Now, we obtain an improved Chen-Ricci inequality in the following:

\begin{th}
\label{th-T-Ricci-imp} Let $(M,g)$ be a Riemannian manifold of dimension $n$
$(n\geq 2)$, $(B,g_{B})$ a Riemannian vector bundle of dimension $(n+s)$
over $M$ and $\zeta $ a $B$-valued symmetric $(1,2)$-tensor field. Suppose
that for any orthonormal basis $\{e_{1},\ldots ,e_{n}\}$ of the tangent
space $T_{p}M$ there is an orthonormal basis $\{e_{n+1},\ldots ,e_{2n+s}\}$
of the Riemannian vector bundle $\left( B,g_{B}\right) $ over $M$ at $p$,
such that
\begin{equation}
\left.
\begin{array}{lll}
\zeta _{jk}^{n+i}=\zeta _{ki}^{n+j}=\zeta _{ij}^{n+k}, & \qquad & i,j,k\in
\left\{ 1,\ldots ,n\right\} ,\medskip \\
\zeta _{ij}^{r}=0, & \qquad & r\in \left\{ 2n+1,\ldots ,2n+s\right\} .%
\end{array}%
\right.  \label{eq-zeta-symm}
\end{equation}%
Let $T$ be a curvature-like tensor satisfying the algebraic Gauss equation $(%
\ref{eq-Gauss-alg})$. Then for any unit vector $X\in T_{p}^{1}M$, we have
\begin{equation}
{\rm Ric}_{T}(X)\leq \frac{n-1}{4n}\,\left\Vert {\rm trace}\,\zeta
\right\Vert ^{2}.  \label{eq-T-Ric-ineq-imp}
\end{equation}%
The equality sign holds for any unit tangent vector at $p$ if and only if
either $\zeta =0$ at $p$ or $n=2$ and
\[
\zeta \left( e_{1},e_{1}\right) =3\mu e_{n+1},\quad \zeta \left(
e_{2},e_{2}\right) =\mu e_{n+1},\quad \zeta \left( e_{1},e_{2}\right) =\mu
e_{n+2}
\]%
for some suitable function $\mu $ with respect to some suitable orthonormal
local frame field.
\end{th}

\noindent {\bf Proof.} Take a point $p\in M$ and an orthonormal frame $%
\left\{ e_{1},...,e_{n}\right\} $ in $T_{p}M$ such that $X=e_{1}$. Then from
(\ref{eq-Gauss-alg}) we have
\begin{equation}
{\rm Ric}_{T}(X)=\sum_{\ell =1}^{n}\sum_{j=2}^{n}\left( \zeta _{11}^{n+\ell
}\zeta _{jj}^{n+\ell }-\left( \zeta _{1j}^{n+\ell }\right) ^{2}\right) .
\label{eq-T-Ric-ineq-imp-1}
\end{equation}%
Since
\[
\sum_{\ell =1}^{n}\sum_{j=2}^{n}\left( \zeta _{1j}^{n+\ell }\right) ^{2}\geq
\sum_{j=2}^{n}\left( \zeta _{1j}^{n+1}\right) ^{2}+\sum_{j=2}^{n}\left(
\zeta _{1j}^{n+j}\right) ^{2},
\]%
therefore (\ref{eq-T-Ric-ineq-imp-1}) gives
\begin{equation}
{\rm Ric}_{T}(X)\leq \sum_{\ell =1}^{n}\sum_{j=2}^{n}\zeta _{11}^{n+\ell
}\zeta _{jj}^{n+\ell }-\sum_{j=2}^{n}\left( \zeta _{1j}^{n+1}\right)
^{2}-\sum_{j=2}^{n}\left( \zeta _{1j}^{n+j}\right) ^{2}.
\label{eq-T-Ric-ineq-imp-2}
\end{equation}%
Using (\ref{eq-zeta-symm}) in (\ref{eq-T-Ric-ineq-imp-2}) we get
\begin{equation}
{\rm Ric}_{T}(X)\leq \sum_{\ell =1}^{n}\sum_{j=2}^{n}\zeta _{11}^{n+\ell
}\zeta _{jj}^{n+\ell }-\sum_{j=2}^{n}\left( \zeta _{11}^{n+j}\right)
^{2}-\sum_{j=2}^{n}\left( \zeta _{jj}^{n+1}\right) ^{2}.
\label{eq-T-Ric-ineq-imp-3}
\end{equation}%
Now, suppose that
\[
f_{1}\left( \zeta _{11}^{n+1},\ldots ,\zeta _{nn}^{n+1}\right) =\zeta
_{11}^{n+1}\sum_{j=2}^{n}\zeta _{jj}^{n+1}-\sum_{j=2}^{n}\left( \zeta
_{jj}^{n+1}\right) ^{2},
\]%
\[
f_{\ell }\left( \zeta _{11}^{n+\ell },\ldots ,\zeta _{nn}^{n+\ell }\right)
=\zeta _{11}^{n+\ell }\sum_{j=2}^{n}\zeta _{jj}^{n+\ell }-\left( \zeta
_{11}^{n+\ell }\right) ^{2},\qquad \ell \in \left\{ 2,\ldots ,n\right\} .
\]%
Since the first component of $\,{\rm trace}\,\zeta \,$ is
\[
\left( {\rm trace}\,\zeta \right) ^{1}=\zeta _{11}^{n+1}+\cdots +\zeta
_{nn}^{n+1},
\]%
by using Lemma~\ref{lem-Deng-1}, we have
\begin{equation}
f_{1}\left( \zeta _{11}^{n+1},\ldots ,\zeta _{nn}^{n+1}\right) \leq \frac{n-1%
}{4n}\,\left( \left( {\rm trace}\,\zeta \right) ^{1}\right) ^{2}.
\label{eq-zeta-f1}
\end{equation}%
Similarly, by Lemma~\ref{lem-Deng-2}, for $2\leq \ell \leq n$, in view of $%
n\geq 2$, we have
\begin{equation}
f_{\ell }\left( \zeta _{11}^{n+\ell },\ldots ,\zeta _{nn}^{n+\ell }\right)
\leq \frac{1}{8}\,\left( \left( {\rm trace}\,\zeta \right) ^{\ell }\right)
^{2}\leq \frac{n-1}{4n}\,\left( \left( {\rm trace}\,\zeta \right) ^{\ell
}\right) ^{2}.  \label{eq-zeta-fr}
\end{equation}%
Now, in view of (\ref{eq-T-Ric-ineq-imp-2}), (\ref{eq-zeta-f1}) and (\ref%
{eq-zeta-fr}) we get
\[
{\rm Ric}_{T}(X)\leq \frac{n-1}{4n}\sum_{\ell =1}^{n}\left( \left( {\rm trace%
}\,\zeta \right) ^{\ell }\right) ^{2}=\frac{n-1}{4n}\,\left\Vert {\rm trace}%
\,\zeta \right\Vert ^{2},
\]%
which gives (\ref{eq-T-Ric-ineq-imp}). \medskip

Now we assume that $n\geq 2$ and the equality sign of (\ref%
{eq-T-Ric-ineq-imp}) is true for all unit vectors $X\in T_{p}^{1}M$. From (%
\ref{eq-zeta-fr}), it follows that $\left( {\rm trace}\,\zeta \right) ^{\ell
}=0$ for $\ell \geq 2$ (or simply choose $e_{n+1}$ parallel to ${\rm trace}%
\,\zeta $). Combining this and Lemma~\ref{lem-Deng-2} we have
\[
\zeta _{1j}^{n+1}=\zeta _{11}^{n+j}=\frac{\left( {\rm trace}\,\zeta \right)
^{j}}{4}=0,\qquad j\geq 2.
\]%
From (\ref{eq-T-Ric-ineq-imp-2}), we get
\[
\zeta _{jk}^{n+1}=0,\qquad j,k\geq 2,\;j\neq k.
\]%
From Lemma~\ref{lem-Deng-1}, the matrix $\left( \zeta _{jk}^{n+1}\right) $
must be diagonal with
\[
\zeta _{11}^{n+1}=\left( n+1\right) \frac{\left( {\rm trace}\,\zeta \right)
^{1}}{2n},\qquad \zeta _{jj}^{n+1}=\frac{\left( {\rm trace}\,\zeta \right)
^{1}}{2n},\qquad j\geq 2.
\]%
Now if we compute ${\rm Ric}_{T}(e_{2})$ as we do for ${\rm Ric}_{T}(X)={\rm %
Ric}_{T}(e_{1})$ in (\ref{eq-T-Ric-ineq-imp-2}), from the equality we get
\[
\zeta _{2j}^{n+\ell }=\zeta _{j\ell }^{n+2}=0,\qquad \ell \neq 2,\;j\neq
2,\;\ell \neq j.
\]%
From the equality and Lemma~\ref{lem-Deng-1}, we obtain
\[
\frac{1}{n+1}\,\zeta _{11}^{n+2}=\zeta _{22}^{n+2}=\cdots =\zeta
_{nn}^{n+2}= \frac{\left( {\rm trace}\,\zeta \right) ^{2}}{2n}=0.
\]%
Since the equality holds for all unit tangent vectors, the argument is also
true for matrices $\left( \zeta _{jk}^{n+\ell }\right) $. Thus, finally we
have
\[
\zeta _{2\ell }^{n+2}=\zeta _{22}^{n+\ell }=\frac{\left( {\rm trace}\,\zeta
\right) ^{\ell }}{2n}=0,\qquad \ell \geq 3.
\]%
Therefore the matrix $\left( \zeta _{jk}^{n+2}\right) $ has only two
possible nonzero entries, that is,
\[
\zeta _{12}^{n+2}=\zeta _{21}^{n+2}=\zeta _{22}^{n+1}=\frac{\left( {\rm trace%
}\,\zeta \right) ^{1}}{2n}.
\]%
Similarly the matrix $\left( \zeta _{jk}^{n+\ell }\right) $ has only two
possible nonzero entries
\[
\zeta _{1\ell }^{n+\ell }=\zeta _{\ell 1}^{n+\ell }=\zeta _{\ell \ell
}^{n+1}=\frac{\left( {\rm trace}\,\zeta \right) ^{1}}{2n},\qquad \ell \geq
3.
\]%
Now, we compute ${\rm Ric}_{T}(e_{2})$ as follows. From (\ref{eq-Gauss-alg}%
), we get
\[
T\left( e_{j},e_{2},e_{2},e_{j}\right) =g_{B}\left( \zeta \left(
e_{j},e_{j}\right) ,\zeta \left( e_{2},e_{2}\right) \right) -g_{B}\left(
\zeta \left( e_{2},e_{j}\right) ,\zeta \left( e_{j},e_{2}\right) \right)
\]%
so we have%
\begin{equation}
T\left( e_{j},e_{2},e_{2},e_{j}\right) =\left( \frac{\left( {\rm trace}%
\,\zeta \right) ^{1}}{2n}\right) ^{2},\qquad j\geq 3.
\label{eq-T(e2,ej,ej,e2)}
\end{equation}%
From
\[
T\left( e_{1},e_{2},e_{2},e_{1}\right) =g_{B}\left( \zeta \left(
e_{1},e_{1}\right) ,\zeta \left( e_{2},e_{2}\right) \right) -g_{B}\left(
\zeta \left( e_{2},e_{1}\right) ,\zeta \left( e_{1},e_{2}\right) \right)
\]%
we get%
\begin{equation}
T\left( e_{1},e_{2},e_{2},e_{1}\right) =\left( n+1\right) \left( \frac{%
\left( {\rm trace}\,\zeta \right) ^{1}}{2n}\right) ^{2}-\left( \frac{\left(
{\rm trace}\,\zeta \right) ^{1}}{2n}\right) ^{2}.  \label{eq-T(e2,e1,e1,e2)}
\end{equation}%
By combining (\ref{eq-T(e2,ej,ej,e2)}) and (\ref{eq-T(e2,e1,e1,e2)}), we get
\[
{\rm Ric}_{T}(e_{2})=\left( n+1\right) \left( \frac{\left( {\rm trace}%
\,\zeta \right) ^{1}}{2n}\right) ^{2}-\left( \frac{\left( {\rm trace}\,\zeta
\right) ^{1}}{2n}\right) ^{2}+\left( n-2\right) \left( \frac{\left( {\rm %
trace}\,\zeta \right) ^{1}}{2n}\right) ^{2},
\]%
which gives
\begin{equation}
{\rm Ric}_{T}(e_{2})=\frac{n-1}{2n^{2}}\left( \left( {\rm trace}\,\zeta
\right) ^{1}\right) ^{2}.  \label{eq-T-Ric(e2)-1}
\end{equation}%
On the other hand from the equality assumption, we have%
\begin{equation}
{\rm Ric}_{T}(e_{2})=\frac{n-1}{4n}\,\left\Vert {\rm trace}\,\zeta
\right\Vert ^{2}=\frac{n-1}{4n}\,\left( \left( {\rm trace}\,\zeta \right)
^{1}\right) ^{2}.  \label{eq-T-Ric(e2)-2}
\end{equation}%
From (\ref{eq-T-Ric(e2)-1}) and (\ref{eq-T-Ric(e2)-2}), it is clear that
\[
\frac{1}{4n^{2}}\,\left( n-1\right) \left( n-2\right) \left( \left( {\rm %
trace}\,\zeta \right) ^{1}\right) ^{2}=0.
\]%
Since $n\neq 1$, we have either $\left( {\rm trace}\,\zeta \right) ^{1}=0$
or $n=2$. If $\left( {\rm trace}\,\zeta \right) ^{1}=0$, then all $\zeta
_{jk}^{n+\ell }$ are zero and hence $\zeta =0$. If $n=2$, then we have
\[
\zeta \left( e_{1},e_{1}\right) =\lambda e_{n+1},\quad \quad \zeta \left(
e_{2},e_{2}\right) =\mu e_{n+1},\quad \quad \zeta \left( e_{1},e_{2}\right)
=\mu e_{n+2},
\]%
with
\[
\lambda =3\mu =\frac{3\left( {\rm trace}\,\zeta \right) ^{1}}{2n}.
\]

The converse is easy to prove by simple computation. $\blacksquare $ \medskip

\section{Lagrangian submanifolds\label{sect-Lagrangian-submanifolds}}

Let $M$ be an $n$-dimensional submanifold of a Riemannian manifold $(%
\widetilde{M},g)$. Then the second fundamental form $\sigma $ of the
immersion is related to the shape operator $A$ by
\[
g\left( \sigma \left( X,Y\right) ,N\right) =g\left( A_{N}X,Y\right) ,
\]%
and the equation of Gauss is given by
\begin{eqnarray}
R(X,Y,Z,W) &=&\widetilde{R}(X,Y,Z,W)+\ g\left( \sigma (X,W),\sigma
(Y,Z)\right)  \label{eq-Gauss} \\
&&-\ g\left( \sigma (X,Z),\sigma (Y,W)\right)  \nonumber
\end{eqnarray}%
for all $X,Y,Z,W\in TM$, where $\widetilde{R}$ and $R$ are the curvature
tensors of $\widetilde{M}$ and $M$ respectively.

\medskip A point $p\in M$ is called a geodesic point if the second
fundamental form $\sigma $ vanishes at $p$. The submanifold is said to be
totally geodesic if every point of $M$ is a geodesic point. A Riemannian
submanifold $M$ is a totally geodesic submanifold of $\widetilde{M}$ if and
only if every geodesic of $M$ is a geodesic of $\widetilde{M}$. The
submanifold $M$ is minimal if the mean curvature vector $H=\frac{1}{n}\,{\rm %
trace}(\sigma )$ vanishes identically. A point $p\in M$ is called an
umbilical point if $\sigma =g\otimes H$ at $p$, that is, the shape operator $%
A_{N}$ is proportional to the identity transformation for all $N\in
T_{p}^{\perp }M$. The submanifold is said to be totally umbilical if every
point of the submanifold is an umbilical point. If the shape operator $A_{H}$
at the mean curvature vector $H$ satisfies $A_{H}X=g(H,H)X$ for every $X\in
TM$, then $M$ is called pseudo-umbilical. Totally umbilical submanifolds are
the simplest submanifolds, which are pseudo-umbilical. Thus for a totally
umbilical submanifold the shape operator $A_{H}$ at $H$ has exactly one
eigenvalue $g(H,H)$; moreover, $A_{N}=0$ for each normal vector $N$
orthogonal to $H$. \medskip

Let $(\widetilde{M},J,g)$ be a $2m$-dimensional almost Hermitian manifold.
If $\widetilde{M}$ is a Kaehler manifold with constant holomorphic sectional
curvature $c$, then it is called a complex space form, denoted by $%
\widetilde{M}(c)$. In this case, the almost complex structure $J$ is
parallel, and the Riemann curvature tensor $\widetilde{R}$ is given by
\begin{eqnarray*}
\widetilde{R}\left( X,Y\right) Z &=& \frac{c}{4}\, \left( g(X,Z)Y - g(Y,Z)X
\right. \\
&&\left. +g(JX,Z)JY-g(JY,Z)JX+2g(JX,Y)JZ\right)
\end{eqnarray*}%
for all vector fields $X,Y,Z$ on $\widetilde{M}$. The model spaces for
complex space forms are the complex Euclidean spaces ${\Bbb C}^{n}$ ($c=0$),
the complex projective spaces ${\Bbb C}P^{n}$ ($c>0$) and the complex
hyperbolic spaces ${\Bbb C}H^{n}$ ($c<0$). \medskip

An $n$-dimensional submanifold $M$ of $(\widetilde{M},J,g)$ is called a
Lagrangian submanifold if the almost complex structure $J$ of $\widetilde{M}$
carries each tangent space of $M$ onto its corresponding normal space, that
is, $J(T_{p}M)=T_{p}^{\perp }M$ for every $p\in M$. \medskip

It is well known from the work of Cartan \cite[p.~231]{Cartan-46-book} that
an $n$-dimensional totally umbilical submanifold of a Euclidean $m$-space is
always an open portion of either an $n$-plane or an $n$-sphere. Totally
umbilical submanifolds, if they exist, are the simplest submanifolds next to
totally geodesic submanifolds in a Riemannian manifold from extrinsic point
of views. However, from Theorem~1 of \cite{Chen-Ogiue-74-MMJ}, it follows
that there exist no totally umbilical Lagrangian submanifolds, of dimension $%
n\geq 2$, in a complex space form except the totally geodesic ones.

\medskip Because of nonexistence of totally umbilical Lagrangian
submanifolds, B.-Y. Chen \cite{Chen-97-Israel} introduced the concept of $H$%
-umbilical Lagrangian submanifolds, which are the simplest Lagrangian
submanifolds next to the totally geodesic ones in complex space forms $%
\widetilde{M}(c)$. By an $H$-umbilical Lagrangian submanifold of a Kaehler
manifold $\widetilde{M}$ we mean a Lagrangian submanifold whose second
fundamental form $\sigma $ assumes the following simple form:
\begin{equation}
\sigma (e_{1},e_{1})=\lambda Je_{1},\quad \sigma (e_{2},e_{2})=\cdots
=\sigma (e_{n},e_{n})=\mu Je_{1},  \label{eq-H-umb-1}
\end{equation}%
\begin{equation}
\sigma (e_{1},e_{j})=\mu Je_{j},\quad \sigma (e_{j},e_{k})=0,\qquad j\neq
k,\quad j,k=2,\ldots ,n  \label{eq-H-umb-2}
\end{equation}%
for some suitable functions $\lambda $ and $\mu $ with respect to some
suitable orthonormal local frame field. \medskip

Now, we apply Theorem~\ref{th-T-Ricci-imp} and obtain an improved Chen-Ricci
inequality for Lagrangian submanifolds of a complex space form in the
following:

\begin{th}
\label{th-T-Ricci-imp-Lag} {\rm (\cite[Theorem~3.1]{Deng-09-IEJG})} Let $M$
be a Lagrangian submanifold of real dimension $n$ $(n\geq 2)$ in a complex
space form $\widetilde{M}(c)$ and $X$ be a unit tangent vector in $%
T_{p}^{1}M $. Then we have%
\begin{equation}
{\rm Ric}(X)\leq \frac{n-1}{4}\,(c+n\left\Vert H(p)\right\Vert ^{2}),
\label{eq-Ricci-imp-Lag}
\end{equation}%
where $H$ is the mean curvature vector of $M$ in $\widetilde{M}(c)$ and $%
{\rm Ric}(X)$ is the Ricci curvature of $M$ at $X$. The equality sign holds
for any unit tangent vector at $p$ if and only if either $p$ is a geodesic
point or $n=2$ and $p$ is an $H$-umbilical point with $\lambda =3\mu $, that
is
\[
\sigma \left( e_{1},e_{1}\right) =3\mu Je_{1},\quad \sigma \left(
e_{2},e_{2}\right) =\mu Je_{1},\quad \sigma \left( e_{1},e_{2}\right) =\mu
Je_{2}
\]%
\noindent for some suitable function $\mu $.
\end{th}

\noindent {\bf Proof.} In (\ref{eq-Gauss-alg}), we set
\[
T\left( X,Y,Z,W\right) =R(X,Y,Z,W)+\frac{c}{4}\,(g(Y,Z)g(X,W)-g(X,Z)g(Y,W))
\]%
with $R$ the Riemannian curvature tensor on $M$, and $\zeta =\sigma $ with $%
\sigma $ the second fundamental form of the immersion of $M$ into $%
\widetilde{M}(c)$. Then we see that
\[
{\rm Ric}_{T}(X)={\rm Ric}(X)-\frac{1}{4}\,(n-1)c
\]%
for all unit vectors $X$. Using this in (\ref{eq-T-Ric-ineq-imp}), we can
complete the proof. $\blacksquare $ \medskip

From Theorem~\ref{th-T-Ricci-imp-Lag}, we have the following

\begin{cor}
{\rm (\cite[Corollary~3.2]{Deng-09-IEJG})} Let $M$ be an $n$-dimensional $%
(n\geq 2)$ Lagrangian submanifold of a complex space form $\widetilde{M}(c)$%
. If
\[
{\rm Ric}(X)=\frac{n-1}{4}\,(c+n\left\Vert H\right\Vert ^{2})
\]%
for all unit tangent vector $X$ of $M$, then either $M$ is a totally
geodesic submanifold in $\widetilde{M}(c)$ or $n=2$ and $M$ is a Lagrangian $%
H$-umbilical surface of $\widetilde{M}(c)$ with $\lambda =3\mu $.
\end{cor}

\begin{ex-new}
(\cite[Example~3.1]{Deng-09-IEJG}) The Whitney $2$-sphere in ${\Bbb C}^{2}$
satisfies the improved Chen-Ricci equality.
\end{ex-new}

\begin{rem-new}
In \cite[Theorem 3.2]{Oprea-05-arXiv}, Oprea proved the improved Chen-Ricci
inequality $(\ref{eq-Ricci-imp-Lag})$ for Lagrangian submanifolds of complex
space forms using optimization techniques on Riemannian submanifolds. Later,
Deng \cite[Theorem~3.1]{Deng-09-IEJG} proved the improved Chen-Ricci
inequality $(\ref{eq-Ricci-imp-Lag})$ by algebraic techniques.
\end{rem-new}

\begin{rem-new}
\label{rem-3.3} $H$-umbilical Lagrangian submanifolds in complex space forms
satisfying the condition $\lambda =3\mu $ have been classified completely
\cite{Chen-97-Israel}. For details, we refer to \cite{Borrelli-CM-95}, \cite%
{Castro-Urbano-93}, \cite{Castro-Urbano-95}, \cite{Chen-96-Edinburg}, \cite[%
p.~331-332]{Chen-00-handbook}, \cite{Chen-08-Romania-Acad} and \cite%
{Chen-Vran-96-MPCPS}.
\end{rem-new}

\begin{prob-new}
To extend Theorem~\ref{th-T-Ricci-imp-Lag} to obtain an improved Chen-Ricci
inequality for Lagrangian submanifolds in Quaternion projective spaces,
which will improve Theorem~3.1 of \cite{Liu-02-AMB}.
\end{prob-new}

\section{Kaehlerian slant submanifolds\label%
{sect-Kaehlerian-slant-submanifolds}}

Let $M$ be a submanifold of an almost Hermitian manifold $(\widetilde{M}%
,J,g) $. We write
\[
JX=PX+FX,\qquad X\in TM,
\]%
where $PX$ and $FX$ are the tangential and the normal components of $JX$,
respectively. Then, $P$ is an endomorphism of the tangent bundle $TM$ and $F$
is a normal bundle valued $1$-form on $TM$. For any nonzero vector $X$
tangent to $M$ at a point $p\in M$, the Wirtinger angle of $X$, denoted by $%
\theta (X)$ is the angle between $JX$ and the tangent space $T_{p}M$. The
submanifold $M$ is called a slant submanifold if $\theta (X)$ is independent
of the choice of $p\in M$ and of $X\in T_{p}M$. The Wirtinger angle of a
slant submanifold is called the slant angle of the slant submanifold. For
slant submanifolds, $P^{2} = t I$, for some $t \in \left[ -1,0\right] $,
where $I$ is the identity transformation of $TM$. Moreover, if $M$ is a
slant submanifold and $\theta $ is the slant angle of $M$, then $t =
-\cos^{2}\theta $. Hence, for a slant submanifold, we have
\[
g\left( PX,PY\right) =\cos ^{2}\theta\, g(X,Y)
\]%
\[
g\left( FX,FY\right) =\sin ^{2}\theta\, g(X,Y)
\]%
for $X, Y$ tangent to $M$. \medskip

We note that a slant submanifold $M$ is $J$-invariant, anti-$J$-invariant,
non-invariant slant or proper slant according as $\theta =0$ $\left(
t=-1\right) $, $\theta =\pi /2$ $\left( t=0\right) $, $\theta \neq 0$ $%
\left( t\neq -1\right) $ or $0\neq \theta \neq \pi /2$ $\left( -1<t=-\cos
^{2}\theta <0\right) $, respectively. \medskip

A proper slant submanifold is said to be Kaehlerian slant if the
endomorphism $P$ is parallel. A Kaehlerian slant submanifold is a Kaehler
manifold with respect to the induced metric and the almost complex structure
$J^{\prime }=\left( \sec \theta \right) J$, where $\theta $ is the slant
angle. Examples of proper slant submanifolds and Kaehlerian slant
submanifolds are given in \cite{Chen-90-slant}. \medskip

For Kaehlerian slant submanifolds in $2n$-dimensional complex space form $%
\widetilde{M}(c)$ we prove the following improved Chen-Ricci inequality.

\begin{th}
\label{th-Ricci-imp-Kaeh-slant} Let $M$ be an $n$-dimensional Kaehlerian
slant submanifold of a $2n$-dimensional complex space form $\widetilde{M}(c)$%
, and $X$ a unit tangent vector in $T_{p}^{1}M$, $p\in M$. Then
\begin{equation}
{\rm Ric}(X)\leq \frac{1}{4}\,\left( \left( n-1\right) n\left\Vert
H\right\Vert ^{2}+(n-1)c+3c\cos ^{2}\theta \right) ,
\label{eq-Ricci-imp-Kaeh-slant}
\end{equation}%
where $H$ is the mean curvature vector of $M$ in $\widetilde{M}(c)$ and $%
{\rm Ric}(X)$ is the Ricci curvature of $M$ at $X$. The equality sign holds
for any unit tangent vector at $p$ if and only if either

\begin{enumerate}
\item[{\rm (a)}] $p$ is a geodesic point or

\item[{\rm (b)}] $n=2$ and
\begin{equation}
\left\{
\begin{array}{c}
\displaystyle \sigma \left( e_{1},e_{1}\right) =3\mu \frac{Fe_{1}} {%
\left\Vert Fe_{1}\right\Vert }=3\mu \csc \theta Fe_{1}, \medskip \\
\displaystyle \sigma \left( e_{2},e_{2}\right) =\mu \frac{Fe_{1}}{\left\Vert
Fe_{1}\right\Vert }=\mu \csc \theta Fe_{1}, \medskip \\
\displaystyle \sigma \left( e_{1},e_{2}\right) =\mu \frac{Fe_{2}}{\left\Vert
Fe_{2}\right\Vert }=\mu \csc \theta Fe_{2}%
\end{array}%
\right.  \label{eq-Ricci-imp-Kaeh-slant-2}
\end{equation}%
\noindent for some suitable function $\mu $, where $\sigma $ is the second
fundamental form of the immersion of $M$ into $\widetilde{M}(c)$.
\end{enumerate}
\end{th}

\noindent {\bf Proof.} In (\ref{eq-Gauss-alg}), we set
\begin{eqnarray*}
&&T\left( X,Y,Z,W\right) =R(X,Y,Z,W)+\frac{c}{4}\,\left\{
g(Y,Z)g(X,W)-g(X,Z)g(Y,W)\right. \\
&&\left. +\,g(PX,Z)g\left( PY,W\right) -g(PY,Z)g\left( PX,W\right)
+2g(PX,Y)g\left( PZ,W\right) \right\}
\end{eqnarray*}%
with $R$ the Riemannian curvature tensor on $M$, and $\zeta =\sigma $. Then
it can be shown that
\[
{\rm Ric}_{T}(X)={\rm Ric}(X)-\frac{1}{4}\,(n-1)c-\frac{3}{4}\,c\cos
^{2}\theta
\]%
for all unit vectors $X$. Using this in (\ref{eq-T-Ric-ineq-imp}), we can
complete the proof. $\blacksquare $

\begin{rem-new}
The inequality (\ref{eq-Ricci-imp-Kaeh-slant}) is an improvement of
Chen-Ricci inequality \cite[inequality {\bf (13)} of Theorem~5.2]%
{Hong-Tri-05-IJPAMS-1} or \cite[Inequality~(2.1) of Theorem~2.1]%
{Mat-Mih-Taz-03}.
\end{rem-new}

We recall that totally umbilical submanifolds, if they exist, are the
simplest submanifolds next to totally geodesic submanifolds in a Riemannian
manifold. From Theorem~1 of \cite{Chen-Ogiue-74-MMJ}, it follows that there
do not exist totally umbilical Lagrangian submanifolds, of dimension $n\geq
2 $, in a complex space form except the totally geodesic ones. We improve
this result in the following

\begin{th}
\label{th-tuL-1} If $M$ is a totally umbilical Lagrangian submanifold of a
Kaehler manifold then either $\dim \left( M\right) =1$ or $M$ is totally
geodesic.
\end{th}

\noindent {\bf Proof.} If $\dim (M)>1$, let $X,Y\in T_{p}M$ such that $%
g(X,Y)=0$ and $g(X,X)=1$. Then
\[
g(H,FY)=g(\sigma (X,X),FY)=g(A_{FY}X,X)=g(A_{FX}Y,X)=g(\sigma (X,Y),FX)=0,
\]%
which shows that $H=0$, and consequently $M$ is totally geodesic. $%
\blacksquare $ \medskip

Recently, Sahin \cite[Theorem~3.1]{Sahin-09-RM} proved that every totally
umbilical proper slant submanifold of a Kaehler manifold is totally
geodesic. Combining the result of Sahin with Theorem~\ref{th-tuL-1}, we get
the following

\begin{th}
\label{th-tuL-2} Every $n$-dimensional ($n\geq 2$) totally umbilical
non-invariant slant submanifold of a $2n$-dimensional Kaehler manifold is
totally geodesic.
\end{th}

However, since the shape operator of every proper slant surface (which is
always Kaehlerian slant) and also every Kaehlerian slant submanifold of a
Kaehler manifold must satisfy another condition
\[
A_{FX}Y=A_{FY}X
\]%
for any $X,Y$ tangent to $M$, there do not exist totally umbilical
Kaehlerian slant submanifold in a Kaehlerian manifold. For these reasons,
B.-Y. Chen \cite{Chen-02-Osaka} studied the simplest slant submanifolds
which satisfy the pseudo-umbilical condition $A_{H}X=g(H,H)X$ and $%
A_{FX}Y=A_{FY}X$, and defined such submanifolds to be slant umbilical
submanifolds, or simply slumbilical submanifolds (although slant
pseudo-umbilical submanifold could be a more correct name). In some sense,
slumbilical submanifolds play the role of totally umbilical submanifolds of
Euclidean space in the family of slant submanifolds. An $n$-dimensional
slant submanifold in a Kaehlerian manifold is a slumbilical submanifold with
slant angle $\theta \in \left( 0,\pi /2\right) $ if its second fundamental
form satisfies \cite{Chen-02-Osaka}
\begin{equation}
\sigma (e_{1},e_{1})=\cdots =\sigma (e_{n},e_{n})=\lambda \csc \theta Fe_{1},
\label{eq-slumb-1}
\end{equation}%
\begin{equation}
\sigma (e_{1},e_{j})=\lambda \csc \theta Fe_{j},\quad \sigma
(e_{j},e_{k})=0,\qquad j\neq k,\quad j,k=2,\ldots ,n  \label{eq-slumb-2}
\end{equation}%
for some suitable function $\lambda $ with respect to some orthonormal frame
field $\left\{ e_{1},\ldots ,e_{n}\right\} $. In \cite{Chen-02-Osaka}, Chen
obtained a complete classification of slumbilical submanifolds in complex
space forms. In fact, there exist twelve families of slumbilical
submanifolds in complex space forms with slant angle $\theta \in \left(
0,\pi /2\right) $. Conversely, every slumbilical submanifold in a complex
space form is given by one of these twelve families. \medskip

Now we return to Theorem~\ref{th-Ricci-imp-Kaeh-slant}, and in view of (\ref%
{eq-Ricci-imp-Kaeh-slant-2}) we observe that a Kaehlerian slant surface,
which is not totally geodesic, satisfying the improved Chen-Ricci equality (%
\ref{eq-Ricci-imp-Kaeh-slant}) is different from slumbilical surfaces. A
proper slumbilical surface cannot satisfy the improved Chen-Ricci equality (%
\ref{eq-Ricci-imp-Kaeh-slant}). Thus, we propose the following

\begin{prob-new}
\label{prob-H-slumb-1} A Kaehlerian slant submanifold $M^{n}$ of a complex
space form $\widetilde{M}(c)$ will be called an $H$-slumbilical submanifold
if its second fundamental form $\sigma $ assumes the following simple form:
\begin{equation}
\sigma (e_{1},e_{1})=\lambda \csc \theta Fe_{1},\quad \sigma
(e_{2},e_{2})=\cdots =\sigma (e_{n},e_{n})=\mu \csc \theta Fe_{1},
\label{eq-H-slumb-1}
\end{equation}%
\begin{equation}
\sigma (e_{1},e_{j})=\mu \csc \theta Fe_{j},\quad \sigma
(e_{j},e_{k})=0,\qquad j\neq k,\quad j,k=2,\ldots ,n  \label{eq-H-slumb-2}
\end{equation}%
for some suitable functions $\lambda $ and $\mu $ with respect to some
suitable orthonormal local frame field $\{e_{1},\ldots ,e_{n}\}$. The
problem is to obtain a complete classification of $H$-slumbilical
submanifolds (at least $H$-slumbilical surfaces) in complex space forms.
\end{prob-new}

\section{$C$-totally real submanifolds\label%
{sect-C-totally-real-submanifolds}}

A differentiable $1$-form $\eta $ on a $(2m+1)$-dimensional differentiable
manifold $\widetilde{M}$ is called a contact form if $\eta \wedge (d\eta
)^{m}\neq 0$ everywhere on $\widetilde{M}$, and $\widetilde{M}$ equipped
with a contact form is a contact manifold. Since rank of $d\eta $ is $2m$,
there exists a unique global vector field $\xi $, called the characteristic
vector field, such that
\begin{equation}
\eta (\xi )=1,\qquad {\pounds }_{\xi }\eta =0,  \label{eq-contact-1}
\end{equation}%
where ${\pounds }_{\xi }$ denotes the Lie differentiation by $\xi $.
Moreover, it is well-known that there exist a Riemannian metric $g$ and a $%
(1,1)$-tensor field $\varphi $ such that
\begin{equation}
\varphi \xi =0,\quad \eta \circ \varphi =0,\quad \eta (X) =g\left( X,\xi
\right) ,  \label{eq-contact-2}
\end{equation}%
\begin{equation}
\varphi ^{2}=-I+\eta \otimes \xi ,\quad d\eta \left( X,Y\right) =g\left(
X,\varphi Y\right) ,  \label{eq-contact-3}
\end{equation}%
\begin{equation}
g(X,Y)=g(\varphi X,\varphi Y)+\eta (X)\eta (Y)  \label{eq-contact-4}
\end{equation}%
for $X,Y\in T\widetilde{M}$. The structure $\left( \eta ,\xi ,\varphi
,g\right) $ is called a contact metric structure and the manifold $%
\widetilde{M}$ endowed with such a structure is said to be a contact metric
manifold. \medskip

The contact metric structure $\left( \eta ,\xi ,\varphi ,g\right) $ on $%
\widetilde{M}$ gives rise to a natural almost Hermitian structure on the
product manifold $\widetilde{M}\times {\Bbb R}$. If this structure is
integrable, then $\widetilde{M}$ is said to be a Sasakian manifold. A
Sasakian manifold is characterized by the condition
\begin{equation}
(\widetilde{\nabla}_{X}\varphi )Y=g(X,Y)\xi -\eta (Y)X,\quad X,Y\in T%
\widetilde{M},  \label{eq-contact-5}
\end{equation}%
where $\widetilde{\nabla}$ is Levi-Civita connection. Also, a contact metric
manifold $\widetilde{M}$ is Sasakian if and only if the curvature tensor $%
\widetilde{R}$ satisfies
\begin{equation}
\widetilde{R}(X,Y)\xi =\eta (Y)X-\eta (X)Y,\quad X,Y\in T\widetilde{M}.
\label{eq-contact-6}
\end{equation}

A plane section in $T_{p}\widetilde{M}$ is called a $\varphi $-section if
there exists a vector $X\in T_{p}\widetilde{M}$ orthogonal to $\xi $ such
that $\{X,\varphi X\}$ span the section. The sectional curvature is called $%
\varphi $-sectional curvature. Just as the sectional curvatures of a
Riemannian manifold determine the curvature completely and the holomorphic
sectional curvatures of a Kaehler manifold determine the curvature
completely, on a Sasakian manifold the $\varphi $-sectional curvatures
determine the curvature completely. Moreover on a Sasakian manifold of
dimension $\geq 5$ if at each point the $\varphi $-sectional curvature is
independent of the choice of $\varphi $-section at the point, it is constant
on the manifold and the curvature tensor is given by
\begin{eqnarray}
\widetilde{R}(X,Y)Z & = & \frac{c+3}{4}\,\left\{ g(Y,Z)X-g(X,Z)Y\right\}
\nonumber \\
&&+\,\frac{c-1}{4}\,\left\{ g(X,\varphi Z)\varphi Y-g(Y,\varphi Z)\varphi
X+2g(X,\varphi Y)\varphi Z\right.  \nonumber \\
&&\qquad \qquad \left. +\ \eta (X)\eta (Z)Y-\eta (Y)\eta (Z)X\right.
\nonumber \\
&&\qquad \qquad \left. +\ g(X,Z)\eta (Y)\xi -(Y,Z)\eta (X)\xi \right\}
\label{eq-Sas-sp-form}
\end{eqnarray}%
for all $X,Y,Z\in T\widetilde{M}$. A Sasakian manifold of constant $\varphi $%
-sectional curvature $c$ is called a Sasakian space form $\widetilde{M}(c)$.
\medskip

A well known result of Tanno \cite{Tanno-69-Tohoku-2} is that a complete
simply connected Sasakian manifold of constant $\varphi $-sectional
curvature $c$ is isometric to one of certain model spaces depending on
whether $c>-3$, $c=-3$ or $c<-3$. The model space for $c>-3$ is a sphere
with a $D$-homothetic deformation of the standard structure. For $c=-3$ the
model space is ${\Bbb R}^{2n+1}$ with the contact form $\eta =\frac{1}{2}%
\,(dz-\sum_{i=1}^{n}y^{i}dx^{i})$ together with the metric $ds^{2}=\eta
\otimes \eta +\frac{1}{4}\sum_{i=1}^{n}((dx^{i})^{2}+(dy^{i})^{2})$. For $%
c<-3$ one has a canonically defined contact metric structure on the product $%
B^{n}\times {\Bbb R}$ where $B^{n}$ is a simply connected bounded domain in $%
C^{n}$ with a Kaehler structure of constant negative holomorphic curvature.
In particular, Sasakian space forms exist for all values of $c$. For more
details we refer to \cite{Blair-02-book}. \medskip

A submanifold $M$ in a contact manifold is called a $C$-totally real
submanifold \cite{YKI-76} if every tangent vector of $M$ belongs to the
contact distribution. Thus, a submanifold $M$ in a contact metric manifold
is a $C$-totally real submanifold if $\xi $ is normal to $M$. A submanifold $%
M$ in an almost contact metric manifold is called anti-invariant \cite{YK-76}
if $\varphi \left( TM\right) \subset T^{\perp }M$. If a submanifold $M$ in a
contact metric manifold is normal to the structure vector field $\xi $, then
it is anti-invariant. Thus $C$-totally real submanifolds in a contact metric
manifold are anti-invariant, as they are normal to $\xi $. \medskip

Now, we apply Theorem~\ref{th-T-Ricci-imp} to get an improved Chen-Ricci
inequality for $C$-totally real submanifolds of a Sasakian space form.

\begin{th}
\label{th-T-Ricci-imp-C-tot} Let $M$ be a $C$-totally real submanifold of
real dimension $n$ $(n\geq 2)$ in a Sasakian space form $\widetilde{M}(c)$
of dimension $2n+1$, and $X$ a unit tangent vector in $T_{p}^{1}M$. Then we
have%
\begin{equation}
{\rm Ric}(X)\leq \frac{n-1}{4}\,(c+3+n\left\Vert H\right\Vert ^{2}),
\label{eq-Ricci-imp-C-tot}
\end{equation}%
where $H$ is the mean curvature vector of $M$ in $\widetilde{M}(c)$ and $%
{\rm Ric}(X)$ is the Ricci curvature of $M^{n}$ at $X$. The equality sign
holds for any unit tangent vector at $p$ if and only if either $p$ is a
geodesic point or $n=2$ and
\[
\sigma \left( e_{1},e_{1}\right) =3\mu \varphi e_{1},\quad \sigma \left(
e_{2},e_{2}\right) =\mu \varphi e_{1},\quad \sigma \left( e_{1},e_{2}\right)
=\mu \varphi e_{2}
\]%
\noindent for some suitable function $\mu $, where $\sigma $ is the second
fundamental form of the immersion of $M$ into $\widetilde{M}(c)$.
\end{th}

\noindent {\bf Proof.} In (\ref{eq-Gauss-alg}), we set $\zeta =\sigma $ and
\begin{eqnarray*}
T\left( X,Y,Z,W\right) &=&R(X,Y,Z,W) \\
&&+\frac{c+3}{4}\,(g(Y,Z)g(X,W)-g(X,Z)g(Y,W)).
\end{eqnarray*}%
Then we see that
\[
{\rm Ric}_{T}(X)={\rm Ric}(X)-\frac{1}{4}\,(n-1)(c+3)
\]%
for all unit vectors $X$. Now, the proof follows by using these data in (\ref%
{eq-T-Ric-ineq-imp}). $\blacksquare $

\begin{rem-new}
The improved Chen-Ricci inequality (\ref{eq-Ricci-imp-C-tot}) is an
improvement of Chen-Ricci inequality \cite[Inequality~(2.1) of Theorem~2.1]%
{Mih-02-JAMS}.
\end{rem-new}

\begin{prob-new}
\label{prob-H-umb-C-tot-real} Like the concept of $H$-umbilical Lagrangian
submanifolds \cite{Chen-97-Israel}, we can define an $H$-umbilical $C$%
-totally real submanifold of a Sasakian space form $\widetilde{M}(c)$. By an
$H$-umbilical $C$-totally real submanifold of a Sasakian manifold $%
\widetilde{M}$ we mean a $C$-totally real submanifold whose second
fundamental form $\sigma $ assumes the following simple form:
\begin{equation}
\sigma (e_{1},e_{1})=\lambda \varphi e_{1},\quad \sigma (e_{2},e_{2})=\cdots
=\sigma (e_{n},e_{n})=\mu \varphi e_{1},  \label{eq-H-umb-C-tot-1}
\end{equation}%
\begin{equation}
\sigma (e_{1},e_{j})=\mu \varphi e_{j},\quad \sigma (e_{j},e_{k})=0,\qquad
j\neq k,\quad j,k=2,\ldots ,n  \label{eq-H-umb-C-tot-2}
\end{equation}%
for some suitable functions $\lambda $ and $\mu $ with respect to some
suitable orthonormal local frame field $\left\{ e_{1},\ldots ,e_{n}\right\} $%
. The problem is to obtain a complete classification of $H$-umbilical $C$%
-totally real submanifolds, or at least $H$-umbilical $C$-totally real
surfaces in Sasakian space forms.
\end{prob-new}

\begin{prob-new}
\label{prob-int-submfd-S-sp} To extend Theorem~\ref{th-T-Ricci-imp-C-tot}
for integral submanifolds of $S$-space forms (cf. \cite{Blair-70}, \cite%
{Kim-Dwi-Tri-07}) and to obtain a complete classification of $H$-umbilical
integral submanifolds of $S$-space forms.
\end{prob-new}

\end{document}